\documentclass{amsart}

\usepackage{amssymb}
\usepackage[cp1251]{inputenc}
\usepackage[english]{babel}
\usepackage{amsmath}
\usepackage{color}
\usepackage{graphics}
\usepackage{epsfig}
\usepackage[all,knot]{xy}

\newtheorem{theorem}{Theorem}
\newtheorem{lemma}{Lemma}
\newtheorem{definition}{Definition}

\newtheorem{cor}{Corollary}

\def\proof{{\bf Proof.}}%\nobreak\noindent}
\def\id{{\hbox{\bf Id}}}
\def\im{{\hbox{\bf Im}\;}}

\def\Out{{\hbox{\bf Out}\;}}

\def\Id{{\mathop{\mbox{ Id}}}}

\begin{document}

\title{On loop extensions and cohomology of loops}

\author{Rolando Jimenez and Quitzeh Morales Mel\'endez}

\date{}

\address{Instituto de Matem\'aticas, Representaci\'on Oaxaca,
Universidad Nacional Aut\'onoma de M\'exico,  Le\'on 2, 68000 Oaxaca de Ju\'arez, Oaxaca,
M\'exico} \email{rolando@matcuer.unam.mx}

\address{Universidad Pedag\'ogica Nacional, unidad 201
Camino a la Zanjita S/N, Col. Noche Buena, Santa Cruz Xoxocotl\'an, Oaxaca.
C.P. 71230}\email{qmoralesme@conacyt.mx}

\subjclass[2010] {Primary: 20N05; secondary: 20J06.}

\keywords{loops, left bol lops, loop extensions, cohomology of groups, cohomology of loops.}

\footnote{The authors were partially supported by CONACyT Grant 151338.}

\begin{abstract}
In this paper are defined cohomology-like groups that classify loop
extensions satisfying a given identity in three variables for association identities, and in two variables for the case of commutativity. It is considered a large amount of identities. 
This groups generalize those defined in works of Nishigori \cite{Nishigori_1} and of Johnson and Leedham-Green 
\cite{Johnson_Leedham-Green}. 
It is computed the number of metacyclic 
extensions for trivial action of the 
quotient on the kernel in one particular case for left Bol 
loops and in general for commutative loops.
\end{abstract}

\maketitle

\section{Some history}

In the paper \cite{Eilenberg-Maclane} it is shown
that equivalent classes of group extensions 
\begin{equation}
0\rightarrow K \rightarrow L\rightarrow Q\rightarrow 0,
\end{equation}
where $K$ is a commutative group and with given action
$Q\times K \longrightarrow K$ by group automorphisms, are
classified by the second cohomology group $H^2(Q,K)$.

%by creating an obstruction theory of the defined kernel
%$(K,\theta)$ where $\theta: Q \rightarrow \Out(K) $ is 
%a homomorphism defining an action of the group $Q$ over 
%the group $K$. The restriction of this action to the 
%center of the group defines a $Q$-module structure and 
%it is shown that this structure completely determines 
%the existence and the number of equivalent extensions
%of the given kernel.

%This is done by defining a group structure on the 
%category of all the kernels (or classes of so called
%similar kernels) with the same center and
%the same $Q$-module structure and showing that
%this group is isomorphic to the 3-dimensional 
%cohomology group of the group $Q$ with coefficients
%in the center. Then, an extension exists if, and only 
%if, the cohomology class of the given kernel is vanishes 
%and the number of extensions is the cardinality of 
%the 2-dimensional cohomology group.

Then, in the paper of Nishigori \cite{Nishigori_1} it is 
shown that, without changing this constructions, it is 
possible to modify the definition of the cohomology
groups in order to include left Bol loop extensions. 

In the next work \cite{Nishigori_2} of Nishigori, the new 
definition of the cohomology groups is extended to the case 
when $Q$ is a loop. This can be done because the action 
is given by a homomorphism $\theta: Q \longrightarrow 
\Out(K)$, and $\Out(K)$ is a group. So, the group
acting on $K$ is some factor group of the loop $Q$. 

\section{The action of the quotient loop on the kernel}

%We consider the analogous construction of 
%kernels $(K,\theta)$, where $K$ is a loop, and 
%$\theta: Q \longrightarrow $, and $G$ is a set of 
%self mappings of $K$ such that de fixed point set
%$K^G$ is de center of $K$, i.e. $K^G=Z(K)$, 
%where $Z(K)=C(K)\cap \Nuc(K) $ is the center of $K$, 
%i.e. the intersection of the Moufang center \linebreak
%$C(K)=\{c\in K | \; ck=kc,\; k\in K \}$ 
%and the nucleus 
%$\Nuc(K)= $\linebreak $=\{n\in K | \; (k_1k_2)n=k_1(k_2n),
%\; k_1,k_2 \in K \}$.
%Note that the center $Z(K)$ is an abelian group.

It is needed to discuss the analogue of an action of 
the quotient loop $Q$ on the kernel $K$. As usual, 
we take a normalized section $s:Q \longrightarrow L $,
i.e. $j\circ s = \Id_Q $, $s(e)=e$, where $e$ denotes
the neutral element. Then, consider the equation
\begin{equation}\label{multiplication}
x_sy_s=(xy)_sf(x,y),
\end{equation}where $x_s$ denotes the image of the element
$x\in Q$ under the section $s$. This gives a set of factors 
$f: Q\times Q \longrightarrow K$. Now we compare with
the similar equation defined by changing the representatives:
\begin{equation}
x_sa\cdot y_sb=(xy)_sf(x,y)\cdot \Phi a\Psi b,
\end{equation}where 
$\Psi,\Phi: K\longrightarrow K$
are some inner maps.

In the case of groups the map $\Psi$
is the identity and the map $\Phi$ is conjugation
by the element $x_s\in L $, and for a left Bol loop extension,
this is an inner map that becomes an automorphism
because the kernel $K$ belongs to the associator 
subgroup of the loop $L$. 

In the general case, it is necessary to use some 
inner maps. Indeed, we have
$$
\begin{array}{rl}
x_sa\cdot y_sb&=(xy)_sf(x,y)\cdot \Phi(a,b)=\\
             &=x_sy_s\cdot \Phi(a,b).
\end{array}
$$
So, we need to take the elements $a,b\in K$
to the right using inner maps and use the invariance 
of the normal subloop $K$ under this maps.

Lets introduce some notations. It is customary to 
denote by $L(x,y)=L_{yx}^{-1}L_{y}L_{x}$ left 
inner maps, by $R(x,y)=R_{xy}^{-1}R_{y}R_{x}$
right inner maps and by $T_x=L_x^{-1}R_x$ middle 
inner maps. Denote by $M(x,y)$ the inner map defined 
by the equation 
$$
(xa)y=x(M(x,y)ay),
$$i.e. 
$M(x,y)=T_y^{-1}L(y,x)^{-1}T_{xy}R(x,y)T_x^{-1}$. 

Then
$$
\begin{array}{rl}
x_sa\cdot y_sb &=(x_sa\cdot y_s)\cdot L(x_sa,y_s)b=\\
               &=(x_s\cdot M(x_s,y_s)a y_s)\cdot L(x_sa,y_s)b=\\
               &=(x_s\cdot y_s T_{y_s}M(x_s,y_s)a )\cdot L(x_sa,y_s)b=\\
               &=(x_s y_s \cdot L(y_s,x_s)T_{y_s}M(x_s,y_s)a )
               \cdot L(x_sa,y_s)b=\\
               &=x_s y_s \cdot [M(x_s y_s,L(x_sa,y_s)b)
                 L(y_s,x_s)T_{y_s}M(x_s,y_s)a \cdot L(x_sa,y_s)b].
\end{array}
$$Note that, if $K$ is contained in the nucleus, i.e.
the elements $a,b \in L$ associate with every element
$l\in L$, then this equation takes the well-known form
$$
\begin{array}{rl}
x_sa\cdot y_sb &=x_s y_s T_{y_s}(a) b=\\
               &=(xy)_s f(x,y)T_{y_s}(a) b,
\end{array}
$$where the restriction to $K$ of the middle map $T_{y_s}$
is an automorphism of this group. Then, family of maps
$T_{y_s}$ defines a (right) $Q$-module structure on the 
nucleus $K$, when this group is commutative. Further in 
this paper it will be assumed that the nucleus $K$ is commutative
and that its image associates with all the elements in the loop $L$.

\section{Identities defining a class of loops}

In order to classify different classes of loop extensions, 
we need to formalize the notion of identity in a loop. 
Let us introduce some notations. Consider the map 
$$
\phi_{i,n}: L^n \longrightarrow L^{n-1},
$$
$n>0,\; 0\leqslant i<n $, given by the formula
$$
\phi_{i,n}(a_0, \dots, a_n)
=(a_0, \dots, a_ia_{i+1},\dots, a_n).
$$ To avoid excessive use of indexes, we 
will denote just $\phi_{i}=\phi_{i,n}$.

In this way, we can see an $n$-term product as a 
composition of the form
$$
\phi_{i_1}\circ\cdots\circ\phi_{i_n}(a_0, \dots, a_n),
$$where $ 0\leqslant i_s < s $ and 
$\phi_{i_1}=\phi_0=\phi: L\times L \rightarrow L$ is 
the loop multiplication. 

\begin{definition}
A repetition operator (or diagonal operator) is
some Cartesian product of maps 
$$
r_{i,j}: L^{n-1} \longrightarrow L^{n},$$
$1\leqslant i,j \leqslant n, i\not = j$
$$
r_{i,j}=
\alpha_{i,j}\circ(\Delta \times \id_{Q^{n-2}})=
\alpha_{1,j}\circ \alpha_{1,j}\circ
(\Delta \times \id_{Q^{n-2}})
$$where 
$\Delta:L\longrightarrow L\times L$
is the diagonal map and $\alpha_{1,j}: L^n \longrightarrow L^{n},\;
1\leqslant j \leqslant n$ is a permutation of coordinates.
\end{definition}
 
\begin{definition}
An association law is an 
equation 
$$
\phi_{i_1}\circ\cdots\circ\phi_{i_n}\circ
r =
\phi_{j_1}\circ\cdots\circ\phi_{j_n}\circ
r,
$$where $r: L^{n-i}\longrightarrow L^{n}$ 
is some fixed repetition operator.
\end{definition}

We remark that this constructions do not apply, 
in general, to loops defined by commutation or 
inversion rules.

\begin{definition}
A \textit{canonical one nest} product is one of the form
$$
\phi_{1}\circ\cdots\circ\phi_{n-1}(a_0, \dots, a_n),
$$i.e. a product obtained starting on the right hand side 
by multiplying by the element on the left.
\end{definition}

In this paper we will be dealing with association laws of the form
$$
\phi_{1}\circ\cdots\circ\phi_{n-1}\circ
r =
\phi_{j_1}\circ\cdots\circ\phi_{j_n}\circ
r.
$$ 

Consider the same association law without
repetitions\footnote{It is clear 
that such association law always is trivial in the sense that
it always imply associativity. 
So, in order to define non trivial loops we need always consider 
identities with repetitions. We consider identities without 
repetition to make more clear the 
construction of the differentials.}
\begin{equation}\label{definingidentity}
\phi_{1}\circ\cdots\circ\phi_{n-1} =
\phi_{j_1}\circ\cdots\circ\phi_{j_n}.
\end{equation}
using the equation (\ref{multiplication})
on an element $(a_1, \dots, a_n)\in Q^n$
the product
$$
\phi_{1}\circ\cdots\circ\phi_{n-1}({a_{1}}_s, \dots, {a_{n}}_s)
$$ 
can be taken to the form
$$
\begin{array}{rl}
[\phi_{1}\circ\cdots\circ\phi_{n-1}(a_{1}, \dots, a_{n})]_s
& \cdot 
f(a_{1},\phi_{1}\circ\cdots\circ\phi_{n-2}(a_{2}, \dots, a_{n}))
\cdot \\
\cdot & f(a_{2},\phi_{1}\circ\cdots\circ\phi_{n-3}(a_{3}, \dots, a_{n}))
\cdots f(a_{n-1}, a_{n}).
\end{array}
$$ The right hand product
$$
\phi_{j_1}\circ\cdots\circ\phi_{j_n}
({a_{1}}_s, \dots, {a_{n}}_s)
$$can be taken to the form
$$
[\phi_{j_1}\circ\cdots\circ\phi_{j_n}(a_{1}, \dots, a_{n})]_s
 \cdot z(a_{1}, \dots, a_{n}),
$$where $k(a_{1}, \dots, a_{n})\in K$.  But 
$$
\phi_{1}\circ\cdots\circ\phi_{n-1}(a_{1}, \dots, a_{n})=
\phi_{j_1}\circ\cdots\circ\phi_{j_n}(a_{1}, \dots, a_{n})
$$ because we have an epimorphism $L\longrightarrow Q$, i.e.
$Q$ is a loop belonging to the class of $L$, so the association
law must be true also for products in $Q$.

So,
\begin{equation}
\label{canonical_1}
\begin{array}{c}
f(a_{1},\phi_{1}\circ\cdots\circ\phi_{n-2}(a_{2}, \dots, a_{n}))
\cdot  f(a_{2},\phi_{1}\circ\cdots\circ\phi_{n-3}(a_{3}, \dots, a_{n}))
\cdots f(a_{n-1}, a_{n})=\\
\\
       =k(a_{1}, \dots, a_{n}).
\end{array}
\end{equation}
This identity, written in additive form, gives an operator
$$\delta : C^2(Q,Z)\longrightarrow C^n(Q,Z).$$
The similarly  defined association law with repetitions gives an operator
$\delta : C^2(Q,Z)\longrightarrow C^{n-r}(Q,Z)$, where 
$r$ is the number of repetitions.

So, we need to compute the element 
$k(a_{1}, \dots, a_{n})\in K$. This will be done in the next section.

\subsection{Identities defined by one nested products}

In this paragraph we compute de element $k(a_{1}, \dots, a_{n})\in K$
for the case when the product 
$\phi_{j_1}\circ\cdots\circ\phi_{j_n}
({a_{1}}_s, \dots, {a_{n}}_s)$ has only one nest. Lets clarify 
this notion. 

\begin{definition}
A product is said to be one nested when we have
$$
j_{m-1}={j_m}\qquad \text{ or } \qquad j_{m-1}={j_m-1}
$$ for every $m=n, \dots ,1$. In other words, after pairing the 
two elements $a_{j_n}, a_{j_n+1}$ we consecutively
pair with the resulting element, either to the right, 
in the case 
$j_{m-1}=j_m$, or to the left, 
in the case 
$j_{m-1}=j_m-1$.
\end{definition}

\begin{definition}\label{canonical1nestlaw}
An association law 
$$
\phi_{1}\circ\cdots\circ\phi_{n-1}\circ
r =
\phi_{j_1}\circ\cdots\circ\phi_{j_n}\circ
r.
$$ 
is said to be canonical one nested if both sides of have one nest.
\end{definition}

Assuming in definition \ref{canonical1nestlaw}
that the pairings on the right side of the equation 
begin to the left, we obtain the identity
\begin{align}\label{reducedcanonical1nest}
&\phi_{1}\cdots\phi_{n-1}=
\phi_{1}^{u_q}\phi_{2}
\cdots\phi_{t_1-1}\phi_{t_1}^{u_2}\phi_{t_1+1}\cdots 
\phi_{t_1+\cdots +t_{p-1}}^{u_1}
\phi_{m-t_p+1}\cdots\phi_{m-1},
\end{align}where 
$t_1+\cdots +t_p -1=m,\;m+u_1+\cdots +u_q=n$.

If all the pairings happen to the left, then we obtain 
the canonical one nest product and the defining identity is 
trivial, that is, there is no defining identity. Also, 
without loss of generality, we can assume that the last 
pairing is to the right, that is, we have the map 
$\phi_{1}^{u_q}$ at the left end of
the composition. This is because, if the pairing finishes to the left, then,
in the defining identity (\ref{definingidentity}),
we have
$$
\phi_{1}\circ\cdots\circ\phi_{m-1}\circ\cdots\circ\phi_{n-1} =
\phi_{1}\circ\cdots\circ\phi_{m-1}\circ\phi_{m}^u\cdots\circ\phi_{j_n}
$$for some $m$ and we can cancel\footnote{We remark that, while
this cancellation is possible, one can not cancel at the left
when powers appear, i.e.
$$
\phi_{1}\circ\phi_{2}=
\phi_{1}^{2}
$$ is a non trivial identity (it defines associativity).} and obtain another identity of the 
form.

As long as the pairings happen to the left, we 
obtain a product given by the left hand side of 
the equation (\ref{canonical_1}). If the pairing is 
to the right, then we need to use the commutation 
rule given by the element $a_i$ in the right 
hand side of the pairing, i.e. we have some factor 
of the form
$$
\begin{array}{c}
f(a_{m-t},\phi_{1}\circ\cdots\circ\phi_{t-2}(a_{m-t+1}, \dots, a_{m}))
\cdots f(a_{m-3},\phi_{1}\circ \phi_{2}(a_{m-2},a_{m-1}, a_{m}))\cdot
\\
\cdot       f(a_{m-2},\phi_{1}(a_{m-1}, a_{m}))
       f(a_{m-1}, a_{m})a_{{m+1}_s}=\\
\\
=a_{{m+1}_s}T_{a_{m+1}}\{
f(a_{m-t},\phi_{1}\circ\cdots\circ\phi_{t-2}(a_{m-t+1}, \dots, a_{m}))
\cdots \\
\cdots f(a_{m-3},\phi_{1}\circ \phi_{2}(a_{m-2},a_{m-1}, a_{m}))
      f(a_{m-2},\phi_{1}(a_{m-1}, a_{m}))\cdot
      \\
\cdot  f(a_{m-1}, a_{m})\}=\\
\\
a_{{m+1}_s}T_{a_{m+1}}
f(a_{m-t},\phi_{1}\circ\cdots\circ\phi_{t-2}(a_{m-t+1}, \dots, a_{m}))
\cdots \\
\cdots T_{a_{m+1}}f(a_{m-3},\phi_{1}\circ \phi_{2}(a_{m-2},a_{m-1}, a_{m}))
      T_{a_{m+1}}f(a_{m-2},\phi_{1}(a_{m-1}, a_{m}))\cdot
      \\
\cdot T_{a_{m+1}} f(a_{m-1}, a_{m})
\end{array}
$$ because $T_{a_{m+1}}$ is a homomorphism. Then, using the 
relation (\ref{multiplication}), this gives a term of the 
form 
$$
\begin{array}{c}
f(\phi_{1}\circ\cdots\circ\phi_{t-1}(a_{m-t}, \dots, a_{m})
,a_{m+1})\cdot
\\
\cdot T_{a_{m+1}}
f(a_{m-t},\phi_{1}\circ\cdots\circ\phi_{t-2}(a_{m-t+1}, \dots, a_{m}))
\cdots \\
\cdots T_{a_{m+1}}f(a_{m-3},\phi_{1}\circ \phi_{2}(a_{m-2},a_{m-1}, a_{m}))
      T_{a_{m+1}}f(a_{m-2},\phi_{1}(a_{m-1}, a_{m}))\cdot
      \\
\cdot T_{a_{m+1}} f(a_{m-1}, a_{m}),
\end{array}
$$and we can continue this process as long as we have 
pairings on the right. The next term of the form
$$
\begin{array}{c}
a_{m+2_s}T_{a_{m+2}}f(\phi_{1}\circ\cdots\circ\phi_{t-1}(a_{m-t}, \dots, a_{m})
,a_{m+1})\cdot
\\
\cdot
T_{a_{m+2}}T_{a_{m+1}}
f(a_{m-t},\phi_{1}\circ\cdots\circ\phi_{t-2}(a_{m-t+1}, \dots, a_{m}))
\cdot
\\
\cdots
\\
\cdot  
T_{a_{m+2}}T_{a_{m+1}}f(a_{m-3},\phi_{1}\circ \phi_{2}(a_{m-2},a_{m-1}, a_{m}))
\cdot
\\
\cdot      
T_{a_{m+2}}T_{a_{m+1}}f(a_{m-2},\phi_{1}(a_{m-1}, a_{m}))\cdot
      \\
\cdot T_{a_{m+2}}T_{a_{m+1}} f(a_{m-1}, a_{m})
\end{array}
$$by (\ref{multiplication}) gives the term
$$
\begin{array}{c}
f(\phi_{1}\phi_{1}\circ\cdots\circ\phi_{t-1}(a_{m-t}, \dots, a_{m},a_{m+1})
,a_{m+2})\cdot
\\
\cdot T_{a_{m+2}}f(\phi_{1}\circ\cdots\circ\phi_{t-1}(a_{m-t}, \dots, a_{m})
,a_{m+1})\cdot
\\
\cdot
T_{a_{m+2}}T_{a_{m+1}}
f(a_{m-t},\phi_{1}\circ\cdots\circ\phi_{t-2}(a_{m-t+1}, \dots, a_{m}))
\cdot
\\
\cdot  
T_{a_{m+2}}T_{a_{m+1}}f(a_{m-t+1},\phi_{1}\circ\cdots\circ\phi_{t-3}(a_{m-t+2}, \dots, a_{m}))
\\
\cdots 
T_{a_{m+2}}T_{a_{m+1}}f(a_{m-1}, a_{m})=
\\
\\
=f(\phi_{1}^2\circ\phi_{2}\circ\cdots\circ\phi_{t-1}(a_{m-t}, \dots, a_{m},a_{m+1})
,a_{m+2})\cdot
\\
\cdot T_{a_{m+2}}f(\phi_{1}\circ\cdots\circ\phi_{t-2}(a_{m-t+1}, \dots, a_{m})
,a_{m+1})\cdot
\\
\cdot
T_{a_{m+2}}T_{a_{m+1}}
f(a_{m-t},\phi_{1}\circ\cdots\circ\phi_{t-2}(a_{m-t+1}, \dots, a_{m}))
\cdot
\\
\cdot  
T_{a_{m+2}}T_{a_{m+1}}f(a_{m-t+1},\phi_{1}\circ\cdots\circ\phi_{t-3}(a_{m-t+2}, \dots, a_{m}))
\\
\cdots 
T_{a_{m+2}}T_{a_{m+1}}f(a_{m-1}, a_{m}).
\end{array}
$$At the end, we obtain

$$
\begin{array}{c}
f(\phi_{1}^{u}\circ\phi_{2}\circ\cdots\circ\phi_{t-1}(a_{m-t}, \dots, a_{m},a_{m+1},\dots,
a_{m+u-1}),a_{m+u})\cdot\\
\cdot
T_{a_{m+u}}f(\phi_{1}^{u-1}\circ\phi_{2}\circ\cdots\circ\phi_{t-1}(a_{m-t}, \dots, a_{m},a_{m+1},\dots,
a_{m+u-2}),a_{m+u-1})
\cdot
\\
\cdot
T_{a_{m+u}}T_{a_{m+u-1}}f(\phi_{1}^{u-2}\circ\phi_{2}\circ\cdots\circ\phi_{t-1}(a_{m-t}, \dots, a_{m},a_{m+1},\dots,
a_{m+u-2}),a_{m+u-3})
\cdots
\\
\cdots 
T_{a_{m+u}}\cdots T_{a_{m+1}}f(a_{m-1}, a_{m}).
\end{array}
$$ Now, in order to pair to the left, we use the 
identity 
\begin{equation}\label{shiftingrule}
\begin{array}{c}
a_{m-t-1}\cdot \phi_{1}^{u}\circ\phi_{2}\circ\cdots\circ\phi_{t-1}(a_{m-t}, \dots, a_{m},a_{m+1},\dots,
a_{m+u-1})
=\\
=\phi_{1}\circ\phi_{2}^{u}\circ\phi_{3}\circ\cdots\circ\phi_{t}(a_{m-t-1},a_{m-t}, \dots, a_{m},a_{m+1},\dots,
a_{m+u-1}).
\end{array}
\end{equation}In this way, from the factor 
$$
\begin{array}{c}
{a_{m-t-1}}_s[\phi_{1}^{u+1}\circ\phi_{2}\circ\cdots\circ\phi_{t-1}(a_{m-t}, \dots, a_{m},a_{m+1},\dots,
a_{m+u-1},a_{m+u})]_s\cdot
\\
\cdot f(\phi_{1}^{u}\circ\phi_{2}\circ\cdots\circ\phi_{t-1}(a_{m-t}, \dots, a_{m},a_{m+1},\dots,
a_{m+u-1}),a_{m+u})\cdot\\
\cdot
T_{a_{m+u}}f(\phi_{1}^{u-1}\circ\phi_{2}\circ\cdots\circ\phi_{t-1}(a_{m-t}, \dots, a_{m},a_{m+1},\dots,
a_{m+u-2}),a_{m+u-1})
\cdot
\\
\cdot
T_{a_{m+u}}T_{a_{m+u-1}}f(\phi_{1}^{u-2}\circ\phi_{2}\circ\cdots\circ\phi_{t-1}(a_{m-t}, \dots, a_{m},a_{m+1},\dots,
a_{m+u-2}),a_{m+u-3})
\cdots
\\
\cdots 
T_{a_{m+u}}\cdots T_{a_{m+1}}f(a_{m-1}, a_{m})
\end{array}
$$by (\ref{multiplication}) we obtain the factor 
$$
\begin{array}{c}
f({a_{m-t-1}},\phi_{1}^{u+1}\circ\phi_{2}\circ\cdots\circ\phi_{t-1}(a_{m-t}, \dots, a_{m},a_{m+1},\dots,
a_{m+u-1},a_{m+u})\cdot
\\
\cdot f(\phi_{1}^{u}\circ\phi_{2}\circ\cdots\circ\phi_{t-1}(a_{m-t}, \dots, a_{m},a_{m+1},\dots,
a_{m+u-1}),a_{m+u})\cdot\\
\cdot
T_{a_{m+u}}f(\phi_{1}^{u-1}\circ\phi_{2}\circ\cdots\circ\phi_{t-1}(a_{m-t}, \dots, a_{m},a_{m+1},\dots,
a_{m+u-2}),a_{m+u-1})
\cdot
\\
\cdot
T_{a_{m+u}}T_{a_{m+u-1}}f(\phi_{1}^{u-2}\circ\phi_{2}\circ\cdots\circ\phi_{t-1}(a_{m-t}, \dots, a_{m},a_{m+1},\dots,
a_{m+u-2}),a_{m+u-3})
\cdots
\\
\cdots 
T_{a_{m+u}}\cdots T_{a_{m+1}}f(a_{m-1}, a_{m}).
\end{array}
$$Further, 
$$
\begin{array}{c}
f({a_{m-t-2}},a_{m-t-1}\phi_{1}^{u+1}\circ\phi_{2}\circ\cdots\circ\phi_{t-1}(a_{m-t}, \dots, a_{m},a_{m+1},\dots,
a_{m+u-1},a_{m+u})\cdot
\\
\cdot
f({a_{m-t-1}},\phi_{1}^{u+1}\circ\phi_{2}\circ\cdots\circ\phi_{t-1}(a_{m-t}, \dots, a_{m},a_{m+1},\dots,
a_{m+u-1},a_{m+u})\cdot
\\
\cdot f(\phi_{1}^{u}\circ\phi_{2}\circ\cdots\circ\phi_{t-1}(a_{m-t}, \dots, a_{m},a_{m+1},\dots,
a_{m+u-1}),a_{m+u})\cdot\\
\cdot
T_{a_{m+u}}f(\phi_{1}^{u-1}\circ\phi_{2}\circ\cdots\circ\phi_{t-1}(a_{m-t}, \dots, a_{m},a_{m+1},\dots,
a_{m+u-2}),a_{m+u-1})
\cdot
\\
\cdot
T_{a_{m+u}}T_{a_{m+u-1}}f(\phi_{1}^{u-2}\circ\phi_{2}\circ\cdots\circ\phi_{t-1}(a_{m-t}, \dots, a_{m},a_{m+1},\dots,
a_{m+u-2}),a_{m+u-3})
\cdots
\\
\cdots 
T_{a_{m+u}}\cdots T_{a_{m+1}}f(a_{m-1}, a_{m})=
\\
=
f({a_{m-t-2}},\phi_{1}\circ\phi_{2}^{u+1}\circ\phi_{3}\circ\cdots\circ\phi_{t-1}(
a_{m-t-1},a_{m-t}, \dots, a_{m},a_{m+1},\dots,
a_{m+u-1},a_{m+u})\cdot
\\
\cdot
f({a_{m-t-1}},\phi_{1}^{u+1}\circ\phi_{2}\circ\cdots\circ\phi_{t-1}(a_{m-t}, \dots, a_{m},a_{m+1},\dots,
a_{m+u-1},a_{m+u})\cdot
\\
\cdot f(\phi_{1}^{u}\circ\phi_{2}\circ\cdots\circ\phi_{t-1}(a_{m-t}, \dots, a_{m},a_{m+1},\dots,
a_{m+u-1}),a_{m+u})\cdot\\
\cdot
T_{a_{m+u}}f(\phi_{1}^{u-1}\circ\phi_{2}\circ\cdots\circ\phi_{t-1}(a_{m-t}, \dots, a_{m},a_{m+1},\dots,
a_{m+u-2}),a_{m+u-1})
\cdot
\\
\cdot
T_{a_{m+u}}T_{a_{m+u-1}}f(\phi_{1}^{u-2}\circ\phi_{2}\circ\cdots\circ\phi_{t-1}(a_{m-t}, \dots, a_{m},a_{m+1},\dots,
a_{m+u-2}),a_{m+u-3})
\cdots
\\
\cdots 
T_{a_{m+u}}\cdots T_{a_{m+1}}f(a_{m-1}, a_{m}).
\end{array}
$$
Then we iterate, that is, we 
substitute $m_2=m-t_1$, $u_1=u$ and repeat this 
computation, till we have finished the pairings
given for indices $u_1,\dots, u_p$ and 
$m_q=m-(t_1+\cdots t_q)$. 

Therefore, he obtain the identity
\begin{equation}
\label{multiplicative_differential1nest}
\begin{array}{c}
f(a_{1},\phi_{1}\cdots\phi_{n-2}(a_{2}, \dots, a_{n}))
\cdot  f(a_{2},\phi_{1}\cdots\phi_{n-3}(a_{3}, \dots, a_{n}))\cdot
\\
\cdots f(a_{n-1}, a_{n})=
\\
= f(\phi_{m-t_1-t_2-\cdots-t_p}^{u_q-1}
    \phi_{m-t_1-t_2-\cdots-t_p+1}
    \cdots
  \phi_{m-t_1-t_2}^{u_2}
   \phi_{m-t_1-t_2+1}\cdots
   \\
   \cdots
   \phi_{m-t_1-1}
   \phi_{m-t_1}^{u_1}\phi_{m-t_1+1}
   \cdots
   \phi_{m-1}\phi_{m}(a_1,\dots ,a_{n-1}),a_{n})
   \\
   T_{a_{n}}
   f(\phi_{m-t_1-t_2-\cdots-t_p}^{u_q-2}
    \phi_{m-t_1-t_2-\cdots-t_p+1}
    \cdots
   \phi_{m-t_1-t_2}^{u_2}
   \phi_{m-t_1-t_2+1}\cdots
   \\
   \cdots
   \phi_{m-t_1-1} 
   \phi_{m-t_1}^{u_1}\phi_{m-t_1+1}
   \cdots
   \phi_{m-1}\phi_{m}(a_1,\dots ,a_{n-2}),a_{n-1})\cdots
   \\
\cdots
   T_{a_{n}}\cdots T_{a_{m+u_1}}f(\phi_{1}^{u_1-1}
   \cdots
   \phi_{m-1}(a_{m-t_1+1}, \dots, a_{m+u_1-2}),a_{m+u_1-1})\cdots
\\
\cdots 
T_{a_{n}}\cdots T_{a_{m+1}}f(a_{m-t_1},\phi_{1}\circ\cdots\circ\phi_{t_1-2}(a_{m-t_1+1}, \dots, a_{m}))
\\
\cdots T_{a_{n}}\cdots T_{a_{m+1}}f(a_{m-2},\phi_{1}(a_{m-1},a_{m}))
\cdot
\\
  \cdot  T_{a_{n}}\cdots T_{a_{m+1}}f(a_{m-1}, a_{m})
\end{array}
\end{equation}(recall that $m-t_1-t_2-\cdots-t_p=1$).

Briefly speaking, left pairings go to the left 
argument of the set of factors and right pairings to the right
argument, on the other argument we have the previous product
(with one less variable) and we need to add $T$ operators for
every missing arguments on the right. More precisely, if we 
have something of the form $f(b,a_{n-l})$ then, in the identity,
the  actual factor is $T_{a_{n}}\cdots T_{a_{n-l+1}}f(b,a_{n-l})$.

By rewriting the previous equation in additive notation, we have the following lemma
\begin{lemma}\label{cocycle_1}
If the identity (\ref{reducedcanonical1nest}) is true in the loop
$L$, then the corresponding set of factors $f\in C^2(Q,K)$ satisfies 
the identity
\begin{equation}
\label{additive_differential1nest}
\begin{array}{c}
f(a_{1},\phi_{1}\cdots\phi_{n-2}(a_{2}, \dots, a_{n}))
+f(a_{2},\phi_{1}\cdots\phi_{n-3}(a_{3}, \dots, a_{n}))+
\\
\cdots +f(a_{n-1}, a_{n})=
\\
= f(\phi_{m-t_1-t_2-\cdots-t_p}^{u_q-1}
    \phi_{m-t_1-t_2-\cdots-t_p+1}
    \cdots
  \phi_{m-t_1-t_2}^{u_2}
   \phi_{m-t_1-t_2+1}\cdots
   \\
   \cdots
   \phi_{m-t_1-1}
   \phi_{m-t_1}^{u_1}\phi_{m-t_1+1}
   \cdots
   \phi_{m-1}\phi_{m}(a_1,\dots ,a_{n-1}),a_{n})+
   \\
+   f(\phi_{m-t_1-t_2-\cdots-t_p}^{u_q-2}
    \phi_{m-t_1-t_2-\cdots-t_p+1}
    \cdots
   \phi_{m-t_1-t_2}^{u_2}
   \phi_{m-t_1-t_2+1}\cdots
   \\
   \cdots
   \phi_{m-t_1-1} 
   \phi_{m-t_1}^{u_1}\phi_{m-t_1+1}
   \cdots
   \phi_{m-1}\phi_{m}(a_1,\dots ,a_{n-2}),a_{n-1})
   a_{n}+\cdots
   \\
   \cdots +f(\phi_{1}^{u_1-1}
   \phi_{m-1}(a_{m-t_1+1}, \dots, a_{m+u_1-2}),
   a_{m+u_1-1})a_{m+u_1}\cdots a_{n}+
\\
\cdots +
f(a_{m-1}, a_{m})a_{m+1}\cdots a_{n} .
\end{array}
\end{equation}
\end{lemma}

This motivates the following definition.

\begin{definition}
The differential $\delta:C^2(Q,Z)\rightarrow C^n(Q,Z)$ associated to 
the identity (\ref{canonical1nestlaw}) is defined by the formula
\begin{equation}
\label{differential1nest}
\begin{array}{c}
\delta f\, (a_1,\dots,a_n)=
\\
= f(\phi_{m-t_1-t_2-\cdots-t_p}^{u_q-1}
    \phi_{m-t_1-t_2-\cdots-t_p+1}
    \cdots
  \phi_{m-t_1-t_2}^{u_2}
   \phi_{m-t_1-t_2+1}\cdots
   \\
   \cdots
   \phi_{m-t_1-1}
   \phi_{m-t_1}^{u_1}\phi_{m-t_1+1}
   \cdots
   \phi_{m-1}\phi_{m}(a_1,\dots ,a_{n-1}),a_{n})+
   \\
+   f(\phi_{m-t_1-t_2-\cdots-t_p}^{u_q-2}
    \phi_{m-t_1-t_2-\cdots-t_p+1}
    \cdots
   \phi_{m-t_1-t_2}^{u_2}
   \phi_{m-t_1-t_2+1}\cdots
   \\
   \cdots
   \phi_{m-t_1-1} 
   \phi_{m-t_1}^{u_1}\phi_{m-t_1+1}
   \cdots
   \phi_{m-1}\phi_{m}(a_1,\dots ,a_{n-2}),a_{n-1})
   a_{n}+\cdots
   \\
   \cdots +f(\phi_{1}^{u_1-1}\cdots
   \phi_{m-1}(a_{m-t_1+1}, \dots, a_{m+u_1-2}),
   a_{m+u_1-1})a_{m+u_1}\cdots a_{n}+
\\
\cdots +
f(a_{m-1}, a_{m})a_{m+1}\cdots a_{n}-
\\
-f(a_{1},\phi_{1}\cdots\phi_{n-2}(a_{2}, \dots, a_{n}))
-f(a_{2},\phi_{1}\cdots\phi_{n-3}(a_{3}, \dots, a_{n}))-
\\
\cdots -f(a_{n-1}, a_{n}).
\end{array}
\end{equation}
\end{definition}

With this definition lemma \ref{cocycle_1} becomes the following.

\begin{lemma}\label{cocycle_2}
If the identity (\ref{reducedcanonical1nest}) is true in the loop
$L$, then the corresponding set of factors $f\in C^2(Q,K)$
is a cocycle with respect to the differential (\ref{differential1nest}).
\end{lemma}

For example, for the case of the left Bol loop identity but without 
repetitions
$$
w(x\cdot yz)=(w\cdot xy)z ,
$$we obtain de identity
$$
f(w,x\cdot yz)f(x,yz)f(y,z)=f(w\cdot xy,z)T_{z_s}f(w,xy)T_{z_s}f(x,y)
$$
or, written in additive form
$$
f(w,x\cdot yz)+f(x,yz)+f(y,z)=f(w\cdot xy,z)+f(w,xy)z+f(x,y)z.
$$This gives the differential
\begin{equation}\label{left-Bol_diff_no_repetitions}
(\delta f)(w,x,y,z)=f(w\cdot xy,z)+f(w,xy)z+f(x,y)z-
f(w,x\cdot yz)-f(x,yz)-f(y,z).
\end{equation}
Now we need to compute the operator $\delta^2$ for the 
cases we have already defined.

Recall that, in dimension one, we have the differential
\begin{equation}\label{differentialdimension1}
(\delta f)(x,y)=f(y)-f(xy)+f(x)y.
\end{equation}
Substituting this in the previous formula we obtain
$$
\begin{array}{l}
(\delta \delta f)(w,x,y,z)=
\\
=f(z)-f((w\cdot xy)z)+f(w\cdot xy)z
\\
+f(xy)z-f(w\cdot xy)z+(f(w)\cdot xy)z
\\
+f(y)z-f(xy)z+f(x)y\cdot z
\\
-f(x\cdot yz)+f(w(x\cdot yz))-f(w)(x\cdot yz)
\\
-f(yz)+f(x\cdot yz)-f(x)\cdot yz
\\
-f(z)+f(yz)-f(y)z.
\end{array}
$$Rearranging, we have 
$$
\begin{array}{l}
(\delta \delta f)(w,x,y,z)=
\\
=f(z)-f(z)\\
-f((w\cdot xy)z)+f(w(x\cdot yz))
\\
+f(w\cdot xy)z-f(w\cdot xy)z
\\
+f(xy)z-f(xy)z
\\
+(f(w)\cdot xy)z-f(w)(x\cdot yz)
\\
+f(y)z-f(y)z
\\
+f(x)y\cdot z-f(x)\cdot yz
\\
-f(x\cdot yz)+f(x\cdot yz)
\\
-f(yz)+f(yz)=
\end{array}
$$
$$
\begin{array}{l}
=-f((w\cdot xy)z)+f(w(x\cdot yz))
\\
+(f(w)\cdot xy)z-f(w)(x\cdot yz)
\\
+f(x)y\cdot z-f(x)\cdot yz.
\end{array}
$$But $(w\cdot xy)z=w(x\cdot yz)$ because
this identity is true over $Q$ by assumption 
and $f(x)y\cdot z=f(x)\cdot yz$, 
$(f(w)\cdot xy)z=f(w)(x\cdot yz)$
because 
$Q$ is acting as group automorphisms, i.e.
the action of $Q$ over $K$ is associative.

Therefore $\delta^2 f\equiv 0$, i.e. $\delta^2= 0$.

\begin{lemma}
The composition $\delta \delta: C^1(Q,K)\rightarrow C^n(Q,K)$
of the differentials (\ref{differential1nest}) and (\ref{differentialdimension1}) equals to zero, i.e. $\delta^2=0$.
\end{lemma}

\begin{equation}
\label{differential_square_1nest}
\begin{array}{c}
\delta \delta f\, (a_1,\dots,a_n)=
\\
= f(a_{n})-
\\
-f(\phi_{m-t_1-t_2-\cdots-t_p}^{u_q-1}
    \phi_{m-t_1-t_2-\cdots-t_p+1}
    \cdots
  \phi_{m-t_1-t_2}^{u_2}
   \phi_{m-t_1-t_2+1}\cdots
   \\
   \cdots
   \phi_{m-t_1-1}
   \phi_{m-t_1}^{u_1}\phi_{m-t_1+1}
   \cdots
   \phi_{m-1}\phi_{m}(a_1,\dots ,a_{n-1}) a_{n})+
\\
+f(\phi_{m-t_1-t_2-\cdots-t_p}^{u_q-1}
    \phi_{m-t_1-t_2-\cdots-t_p+1}
    \cdots
  \phi_{m-t_1-t_2}^{u_2}
   \phi_{m-t_1-t_2+1}\cdots
   \\
   \cdots
   \phi_{m-t_1-1}
   \phi_{m-t_1}^{u_1}\phi_{m-t_1+1}
   \cdots
   \phi_{m-1}\phi_{m}(a_1,\dots ,a_{n-1}))a_{n}+
   \\
   \\
+   f(a_{n-1})a_{n}-
\\
-f(\phi_{m-t_1-t_2-\cdots-t_p}^{u_q-2}
    \phi_{m-t_1-t_2-\cdots-t_p+1}
    \cdots
   \phi_{m-t_1-t_2}^{u_2}
   \phi_{m-t_1-t_2+1}\cdots
   \\
   \cdots
   \phi_{m-t_1-1} 
   \phi_{m-t_1}^{u_1}\phi_{m-t_1+1}
   \cdots
   \phi_{m-1}\phi_{m}(a_1,\dots ,a_{n-2})a_{n-1})a_{n}+
\\
+f(\phi_{m-t_1-t_2-\cdots-t_p}^{u_q-2}
    \phi_{m-t_1-t_2-\cdots-t_p+1}
    \cdots
   \phi_{m-t_1-t_2}^{u_2}
   \phi_{m-t_1-t_2+1}\cdots
   \\
   \cdots
   \phi_{m-t_1-1} 
   \phi_{m-t_1}^{u_1}\phi_{m-t_1+1}
   \cdots
   \phi_{m-1}\phi_{m}(a_1,\dots ,a_{n-2}))a_{n-1}a_{n}+\cdots
   \\
   \\
   \cdots +f(a_{m+u_1-1})a_{m+u_1}\cdots a_{n}-
   \\
   -f(\phi_{1}^{u_1-1}\cdots
   \phi_{m-1}(a_{m-t_1+1}, \dots, a_{m+u_1-2}) a_{m+u_1-1})a_{m+u_1}\cdots a_{n}+
   \\
   +f(\phi_{1}^{u_1-1}\cdots
   \phi_{m-1}(a_{m-t_1+1}, \dots, a_{m+u_1-2}))a_{m+u_1-1}a_{m+u_1}\cdots a_{n}+
\\
\\
\cdots +
f(a_{m})a_{m+1}\cdots a_{n}
-f(a_{m-1}a_{m})a_{m+1}\cdots a_{n}+
f(a_{m-1})a_{m}a_{m+1}\cdots a_{n}-
\\
\\
-f(\phi_{1}\cdots\phi_{n-2}(a_{2}, \dots, a_{n}))
+f(a_{1}\phi_{1}\cdots\phi_{n-2}(a_{2}, \dots, a_{n}))
-f(a_{1})\phi_{1}\cdots\phi_{n-2}(a_{2}, \dots, a_{n})
\\
-f(\phi_{1}\cdots\phi_{n-3}(a_{3}, \dots, a_{n}))
+f(a_{2}\phi_{1}\cdots\phi_{n-3}(a_{3}, \dots, a_{n}))-
f(a_{2})\phi_{1}\cdots\phi_{n-3}(a_{3}, \dots, a_{n})+
\\
\cdots -f( a_{n})+f(a_{n-1} a_{n})-f(a_{n-1})a_{n}.
\end{array}
\end{equation}Now we use the associativity of the actions, i.e.
$ z(ab)=(za)b,\; z\in Z, a,b\in Q$, then
\begin{equation}
\begin{array}{c}
\delta \delta f\, (a_1,\dots,a_n)=
\\
= f(a_{n})-
\\
-f(\phi_{m-t_1-t_2-\cdots-t_p}^{u_q-1}
    \phi_{m-t_1-t_2-\cdots-t_p+1}
    \cdots
  \phi_{m-t_1-t_2}^{u_2}
   \phi_{m-t_1-t_2+1}\cdots
   \\
   \cdots
   \phi_{m-t_1-1}
   \phi_{m-t_1}^{u_1}\phi_{m-t_1+1}
   \cdots
   \phi_{m-1}\phi_{m}(a_1,\dots ,a_{n-1}) a_{n})+
\\
+f(\phi_{m-t_1-t_2-\cdots-t_p}^{u_q-1}
    \phi_{m-t_1-t_2-\cdots-t_p+1}
    \cdots
  \phi_{m-t_1-t_2}^{u_2}
   \phi_{m-t_1-t_2+1}\cdots
   \\
   \cdots
   \phi_{m-t_1-1}
   \phi_{m-t_1}^{u_1}\phi_{m-t_1+1}
   \cdots
   \phi_{m-1}\phi_{m}(a_1,\dots ,a_{n-1}))a_{n}+
   \\
   \\
+   f(a_{n-1})a_{n}-
\\
-f(\phi_{m-t_1-t_2-\cdots-t_p}^{u_q-2}
    \phi_{m-t_1-t_2-\cdots-t_p+1}
    \cdots
   \phi_{m-t_1-t_2}^{u_2}
   \phi_{m-t_1-t_2+1}\cdots
   \\
   \cdots
   \phi_{m-t_1-1} 
   \phi_{m-t_1}^{u_1}\phi_{m-t_1+1}
   \cdots
   \phi_{m-1}\phi_{m}(a_1,\dots ,a_{n-2})a_{n-1})a_{n}+
\\
+f(\phi_{m-t_1-t_2-\cdots-t_p}^{u_q-2}
    \phi_{m-t_1-t_2-\cdots-t_p+1}
    \cdots
   \phi_{m-t_1-t_2}^{u_2}
   \phi_{m-t_1-t_2+1}\cdots
   \\
   \cdots
   \phi_{m-t_1-1} 
   \phi_{m-t_1}^{u_1}\phi_{m-t_1+1}
   \cdots
   \phi_{m-1}\phi_{m}(a_1,\dots ,a_{n-2}))a_{n-1}a_{n}+\cdots
   \\
   \\
   \cdots +f(a_{m+u_1-1})a_{m+u_1}\cdots a_{n}-
   \\
   -f(\phi_{1}^{u_1-1}\cdots
   \phi_{m-1}(a_{m-t_1+1}, \dots, a_{m+u_1-2}) a_{m+u_1-1})a_{m+u_1}\cdots a_{n}+
   \\
   +f(\phi_{1}^{u_1-1}\cdots
   \phi_{m-1}(a_{m-t_1+1}, \dots, a_{m+u_1-2}))a_{m+u_1-1}a_{m+u_1}\cdots a_{n}+
\\
\\
\cdots +
f(a_{m})a_{m+1}\cdots a_{n}
-f(a_{m-1}a_{m})a_{m+1}\cdots a_{n}+
f(a_{m-1})a_{m}a_{m+1}\cdots a_{n}-
\\
\\
-f(\phi_{1}\cdots\phi_{n-2}(a_{2}, \dots, a_{n}))
+f(a_{1}\phi_{1}\cdots\phi_{n-2}(a_{2}, \dots, a_{n}))
-f(a_{1})a_{2} \cdots a_{n}
\\
-f(\phi_{1}\cdots\phi_{n-3}(a_{3}, \dots, a_{n}))
+f(a_{2}\phi_{1}\cdots\phi_{n-3}(a_{3}, \dots, a_{n}))-
f(a_{2})a_{3} \cdots a_{n}+
\\
\cdots -f( a_{n})+f(a_{n-1} a_{n})-f(a_{n-1})a_{n}.
\end{array}
\end{equation}Using the identities
\begin{equation}
\label{leftshiftingrule}
a_{1}\phi_{1}\cdots\phi_{n-2}(a_{2}, \dots, a_{n})=
\phi_{1}\cdots\phi_{n-1}(a_{1},a_{2}, \dots, a_{n})
\end{equation} and
\begin{equation}
\label{rightshiftingrule}
\phi_{1}^{l-1}\phi_{j_{n-2+l}}\cdots\phi_{j_{n-2}}(a_{1}, \dots, a_{n-1})a_{n}=
\phi_{1}^{l}\phi_{j_{n-2+l}}\cdots\phi_{n-2}(a_{1}, \dots, a_{n-1}, a_{n})
\end{equation}
we obtain
\begin{equation}
\begin{array}{c}
\delta \delta f\, (a_1,\dots,a_n)=
\\
= {\color{blue}f(a_{n})}-
\\
-f(\phi_{m-t_1-t_2-\cdots-t_p}^{u_q}
    \phi_{m-t_1-t_2-\cdots-t_p+1}
    \cdots
  \phi_{m-t_1-t_2}^{u_2}
   \phi_{m-t_1-t_2+1}\cdots
   \\
   \cdots
   \phi_{m-t_1-1}
   \phi_{m-t_1}^{u_1}\phi_{m-t_1+1}
   \cdots
   \phi_{m-1}\phi_{m}(a_1,\dots ,a_{n-1}, a_{n}))+
\\
+f(\phi_{m-t_1-t_2-\cdots-t_p}^{u_q-1}
    \phi_{m-t_1-t_2-\cdots-t_p+1}
    \cdots
  \phi_{m-t_1-t_2}^{u_2}
   \phi_{m-t_1-t_2+1}\cdots
   \\
   \cdots
   \phi_{m-t_1-1}
   \phi_{m-t_1}^{u_1}\phi_{m-t_1+1}
   \cdots
   \phi_{m-1}\phi_{m}(a_1,\dots ,a_{n-1}))a_{n}+
   \\
   \\
+   {\color{red}f(a_{n-1})a_{n}}-
\\
-f(\phi_{m-t_1-t_2-\cdots-t_p}^{u_q-1}
    \phi_{m-t_1-t_2-\cdots-t_p+1}
    \cdots
   \phi_{m-t_1-t_2}^{u_2}
   \phi_{m-t_1-t_2+1}\cdots
   \\
   \cdots
   \phi_{m-t_1-1} 
   \phi_{m-t_1}^{u_1}\phi_{m-t_1+1}
   \cdots
   \phi_{m-1}\phi_{m}(a_1,\dots ,a_{n-2},a_{n-1}))a_{n}+
\\
+f(\phi_{m-t_1-t_2-\cdots-t_p}^{u_q-2}
    \phi_{m-t_1-t_2-\cdots-t_p+1}
    \cdots
   \phi_{m-t_1-t_2}^{u_2}
   \phi_{m-t_1-t_2+1}\cdots
   \\
   \cdots
   \phi_{m-t_1-1} 
   \phi_{m-t_1}^{u_1}\phi_{m-t_1+1}
   \cdots
   \phi_{m-1}\phi_{m}(a_1,\dots ,a_{n-2}))a_{n-1}a_{n}+\cdots
   \\
   \\
   \cdots +f(a_{m+u_1-1})a_{m+u_1}\cdots a_{n}-
   \\
   -f(\phi_{1}^{u_1}\cdots
   \phi_{m-1}(a_{m+u_1-1},a_{m-t_1+1}, \dots, a_{m+u_1-2},
   a_{m+u_1-2}))a_{m+u_1}\cdots a_{n}+
   \\
   +f(\phi_{1}^{u_1-1}\cdots
   \phi_{m-1}(a_{m-t_1+1}, \dots, a_{m+u_1-2}))a_{m+u_1-1}a_{m+u_1}\cdots a_{n}+
\\
\\
\cdots +
f(a_{m})a_{m+1}\cdots a_{n}
-f(a_{m-1}a_{m})a_{m+1}\cdots a_{n}+
f(a_{m-1})a_{m}a_{m+1}\cdots a_{n}-
\\
\\
-f(\phi_{1}\cdots\phi_{n-2}(a_{2}, \dots, a_{n}))
+f(\phi_{1}\cdots\phi_{n-1}(a_{1},a_{2}, \dots, a_{n}))
-f(a_{1})a_{2} \cdots a_{n}
\\
-f(\phi_{1}\cdots\phi_{n-3}(a_{3}, \dots, a_{n}))
+f(\phi_{1}\cdots\phi_{n-2}(a_{2},a_{3}, \dots, a_{n}))-
f(a_{2})a_{3} \cdots a_{n}+
\\
\cdots -{\color{blue}f(a_{n})}+f(a_{n-1} a_{n})-{\color{red}f(a_{n-1})a_{n}}.
\end{array}
\end{equation}Note that the sum in the definition (\ref{differential1nest}) of 
$\delta f$ is separated in positive
and negative terms. So, if we preserve this separation
and the 3-term subdivision that comes from the substitution
of the one dimensional differential, then we observe that the
first summand of the first positive term cancels 
with the first summand of the last negative term. 
The first summand of the second positive term cancels 
with the third summand of the last negative term. Also, 
the second summand of the first positive term cancels with
the second sumand of the first positive term, that is 
\begin{align*}
&f(\phi_{m-t_1-t_2-\cdots-t_p}^{u_q}
    \phi_{m-t_1-t_2-\cdots-t_p+1}
    \cdots
  \phi_{m-t_1-t_2}^{u_2}
   \phi_{m-t_1-t_2+1}\cdots
   \\
   &\cdots
   \phi_{m-t_1-1}
   \phi_{m-t_1}^{u_1}\phi_{m-t_1+1}
   \cdots
   \phi_{m-1}\phi_{m}(a_1,\dots ,a_{n-1}, a_{n}))=
   \\
 &f(\phi_{1}\cdots\phi_{n-1}(a_{1},a_{2}, \dots, a_{n}))  
\end{align*}
because the defining identity
\begin{align*}
&\phi_{m-t_1-t_2-\cdots-t_p}^{u_q}
    \phi_{m-t_1-t_2-\cdots-t_p+1}
    \cdots
  \phi_{m-t_1-t_2}^{u_2}
   \phi_{m-t_1-t_2+1}\cdots
   \\
   &\cdots
   \phi_{m-t_1-1}
   \phi_{m-t_1}^{u_1}\phi_{m-t_1+1}
   \cdots
   \phi_{m-1}\phi_{m}(a_1,\dots ,a_{n-1}, a_{n})=
   \\
 &\phi_{1}\cdots\phi_{n-1}(a_{1},a_{2}, \dots, a_{n})  
\end{align*}must be true over $Q$, being a quotient of $L$:
\begin{equation}
\begin{array}{c}
\delta \delta f\, (a_1,\dots,a_n)=
\\
= {\color{blue}f(a_{n})}-
\\
-{\color{magenta}f(\phi_{m-t_1-t_2-\cdots-t_p}^{u_q}
    \phi_{m-t_1-t_2-\cdots-t_p+1}
    \cdots
  \phi_{m-t_1-t_2}^{u_2}
   \phi_{m-t_1-t_2+1}\cdots}
   \\
   {\color{magenta}\cdots
   \phi_{m-t_1-1}
   \phi_{m-t_1}^{u_1}\phi_{m-t_1+1}
   \cdots
   \phi_{m-1}\phi_{m}(a_1,\dots ,a_{n-1}, a_{n}))}+
\\
+f(\phi_{m-t_1-t_2-\cdots-t_p}^{u_q-1}
    \phi_{m-t_1-t_2-\cdots-t_p+1
    \cdots
  \phi_{m-t_1-t_2}^{u_2}
   \phi_{m-t_1-t_2+1}\cdots}
   \\
   \cdots
   \phi_{m-t_1-1}
   \phi_{m-t_1}^{u_1}\phi_{m-t_1+1}
   \cdots
   \phi_{m-1}\phi_{m}(a_1,\dots ,a_{n-1}))a_{n}+
   \\
   \\
+   {\color{red}f(a_{n-1})a_{n}}-
\\
-f(\phi_{m-t_1-t_2-\cdots-t_p}^{u_q-1}
    \phi_{m-t_1-t_2-\cdots-t_p+1}
    \cdots
   \phi_{m-t_1-t_2}^{u_2}
   \phi_{m-t_1-t_2+1}\cdots
   \\
   \cdots
   \phi_{m-t_1-1} 
   \phi_{m-t_1}^{u_1}\phi_{m-t_1+1}
   \cdots
   \phi_{m-1}\phi_{m}(a_1,\dots ,a_{n-2},a_{n-1}))a_{n}+
\\
+f(\phi_{m-t_1-t_2-\cdots-t_p}^{u_q-2}
    \phi_{m-t_1-t_2-\cdots-t_p+1}
    \cdots
   \phi_{m-t_1-t_2}^{u_2}
   \phi_{m-t_1-t_2+1}\cdots
   \\
   \cdots
   \phi_{m-t_1-1} 
   \phi_{m-t_1}^{u_1}\phi_{m-t_1+1}
   \cdots
   \phi_{m-1}\phi_{m}(a_1,\dots ,a_{n-2}))a_{n-1}a_{n}+\cdots
   \\
   \\
   \cdots +f(a_{m+u_1-1})a_{m+u_1}\cdots a_{n}-
   \\
   -f(\phi_{1}^{u_1}\cdots
   \phi_{m-1}(a_{m+u_1-1},a_{m-t_1+1}, \dots, a_{m+u_1-2},
   a_{m+u_1-2}))a_{m+u_1}\cdots a_{n}+
   \\
   +f(\phi_{1}^{u_1-1}\cdots
   \phi_{m-1}(a_{m-t_1+1}, \dots, a_{m+u_1-2}))a_{m+u_1-1}a_{m+u_1}\cdots a_{n}+
\\
\\
\cdots +
f(a_{m})a_{m+1}\cdots a_{n}
-f(a_{m-1}a_{m})a_{m+1}\cdots a_{n}+
f(a_{m-1})a_{m}a_{m+1}\cdots a_{n}-
\\
\\
-f(\phi_{1}\cdots\phi_{n-2}(a_{2}, \dots, a_{n}))
+{\color{magenta}f(\phi_{1}\cdots\phi_{n-1}(a_{1},a_{2}, \dots, a_{n}))}
-f(a_{1})a_{2} \cdots a_{n}
\\
-f(\phi_{1}\cdots\phi_{n-3}(a_{3}, \dots, a_{n}))
+f(\phi_{1}\cdots\phi_{n-2}(a_{2},a_{3}, \dots, a_{n}))-
f(a_{2})a_{3} \cdots a_{n}+
\\
\cdots -{\color{blue}f(a_{n})}+f(a_{n-1} a_{n})-{\color{red}f(a_{n-1})a_{n}}.
\end{array}
\end{equation}
These are cancelations of the positive side with the negative
side of the sum. 

There is also cancelations within each side.
The third summand of the first positive term cancels
with the second summand of the second positive term.
The second sumand of the last negative term cancels with
the first summand of the previous negative term:
\begin{equation}
\begin{array}{c}
\delta \delta f\, (a_1,\dots,a_n)=
\\
= f(a_{n})-
\\
-f(\phi_{m-t_1-t_2-\cdots-t_p}^{u_q}
    \phi_{m-t_1-t_2-\cdots-t_p+1}
    \cdots
  \phi_{m-t_1-t_2}^{u_2}
   \phi_{m-t_1-t_2+1}\cdots
   \\
   \cdots
   \phi_{m-t_1-1}
   \phi_{m-t_1}^{u_1}\phi_{m-t_1+1}
   \cdots
   \phi_{m-1}\phi_{m}(a_1,\dots ,a_{n-1}, a_{n}))+
\\
+{\color{magenta}f(\phi_{m-t_1-t_2-\cdots-t_p}^{u_q-1}
    \phi_{m-t_1-t_2-\cdots-t_p+1}
    \cdots
  \phi_{m-t_1-t_2}^{u_2}
   \phi_{m-t_1-t_2+1}\cdots}
   \\
   {\color{magenta}\cdots
   \phi_{m-t_1-1}
   \phi_{m-t_1}^{u_1}\phi_{m-t_1+1}
   \cdots
   \phi_{m-1}\phi_{m}(a_1,\dots ,a_{n-1}))a_{n}}+
   \\
   \\
+   f(a_{n-1})a_{n}-
\\
-{\color{magenta}f(\phi_{m-t_1-t_2-\cdots-t_p}^{u_q-1}
    \phi_{m-t_1-t_2-\cdots-t_p+1}
    \cdots
   \phi_{m-t_1-t_2}^{u_2}
   \phi_{m-t_1-t_2+1}\cdots}
   \\
   \cdots
   {\color{magenta}\phi_{m-t_1-1} 
   \phi_{m-t_1}^{u_1}\phi_{m-t_1+1}
   \cdots
   \phi_{m-1}\phi_{m}(a_1,\dots ,a_{n-2},a_{n-1}))a_{n}}+
\\
+f(\phi_{m-t_1-t_2-\cdots-t_p}^{u_q-2}
    \phi_{m-t_1-t_2-\cdots-t_p+1}
    \cdots
   \phi_{m-t_1-t_2}^{u_2}
   \phi_{m-t_1-t_2+1}\cdots
   \\
   \cdots
   \phi_{m-t_1-1} 
   \phi_{m-t_1}^{u_1}\phi_{m-t_1+1}
   \cdots
   \phi_{m-1}\phi_{m}(a_1,\dots ,a_{n-2}))a_{n-1}a_{n}+\cdots
   \\
   \\
   \cdots +f(a_{m+u_1-1})a_{m+u_1}\cdots a_{n}-
   \\
   -f(\phi_{1}^{u_1}\cdots
   \phi_{m-1}(a_{m+u_1-1},a_{m-t_1+1}, \dots, a_{m+u_1-2},
   a_{m+u_1-2}))a_{m+u_1}\cdots a_{n}+
   \\
   +f(\phi_{1}^{u_1-1}\cdots
   \phi_{m-1}(a_{m-t_1+1}, \dots, a_{m+u_1-2}))a_{m+u_1-1}a_{m+u_1}\cdots a_{n}+
\\
\\
\cdots +
f(a_{m})a_{m+1}\cdots a_{n}
-f(a_{m-1}a_{m})a_{m+1}\cdots a_{n}+
f(a_{m-1})a_{m}a_{m+1}\cdots a_{n}-
\\
\\
-f(\phi_{1}\cdots\phi_{n-2}(a_{2}, \dots, a_{n}))
+f(\phi_{1}\cdots\phi_{n-1}(a_{1},a_{2}, \dots, a_{n}))
-f(a_{1})a_{2} \cdots a_{n}
\\
-f(\phi_{1}\cdots\phi_{n-3}(a_{3}, \dots, a_{n}))
+f(\phi_{1}\cdots\phi_{n-2}(a_{2},a_{3}, \dots, a_{n}))-
f(a_{2})a_{3} \cdots a_{n}+
\\
\cdots 
-{\color{red}f(\phi_{1}(a_{n-1},a_{n}))}
+f(\phi_{1}\phi_2(a_{n-2},a_{n-1},a_{n}))
-f(a_{n-2})a_{n-1}a_{n}-
\\
-f(a_{n})+{\color{red}f(a_{n-1} a_{n})}-f(a_{n-1})a_{n}.
\end{array}
\end{equation}
Following these instructions (excepting the one involving 
the defining identity, which is used only one time) one 
can show that all terms cancel. Therefore $\delta^2=0$
for every one nested identity. $\hfill \square$

Note that the previous computations do not alter any variables,
we have $\delta^2f=0$ also for a set of factors 
$f: Q^{n-r}\rightarrow Z$ with repetitions.

\subsection{Cocycles defining a non-associative loop}

Notice that, from the identity (\ref{reducedcanonical1nest}) 
it is clear that, after substituting some variables
by the neutral element of the loop, we can obtain the identity
\begin{align*}
&\phi_{1}\circ\phi_{2}=
\phi_{1}^{2},
\end{align*}which defines associativity. That is, all these
identities are trivial when considered without repetitions.
We only worked with this identities to construct the differential
and to show that $\delta^2=0$ in any case.   

In general, it is a complicated problem to show that the identity 
analogous to (\ref{reducedcanonical1nest}) but with repetitions 
defines a non-associative loop. This is a word problem in the 
non-associative context. 

For this reason we will be mostly dealing with identities with three variables.

In this case we have the identity
\begin{equation}\label{reducedcanonical1nest3variables}
\begin{array}{r}
(\phi_{1}\cdots\phi_{n-1})\;\circ r =\\
(\phi_{1}^{u_q}\phi_{2}
\cdots\phi_{t_1-1}\phi_{t_1}^{u_2}\phi_{t_1+1}\cdots\qquad\qquad\qquad\qquad\\
\cdots\phi_{t_1+\cdots +t_{p-1}}^{u_1}
\phi_{m-t_p+1}\cdots\phi_{m-1})\circ r,
\end{array}
\end{equation}where 
$t_1+\cdots +t_p -1=m,\;m+u_1+\cdots +u_q=n\geqslant 4$ and 
$r:Q^3\rightarrow Q^n$ is a repetition operator, that is
$$
r(a_1,a_2,a_3)=\alpha(\;\underbrace{a_1,\dots,a_1}_\text{$i$ times}\;,
                     \;\underbrace{a_2,\dots,a_2}_\text{$j$ times}\;,
                      \underbrace{a_3,\dots,a_3}_\text{$n-i-j$ times})
$$where $\alpha\in S_n$ is some permutation. In this case, every 
elimination of variables, in other words, every substitution of 
a variable by de neutral element $e\in L$ gives an identity in two
variables. Such identity can not define associativity for every
element in the loop $L$ and defines some association rule for pairs of 
elements in this loop. 

\begin{lemma}\label{cartesian_extension}
The product rule
\begin{equation}\label{productrule1nest}
(a,w)\cdot(b,z) = (ab,f(a,b)+ wb+z)
\end{equation}for elements 
$(a,w),(b,z)\in L\times Z$, an action
$\theta_0:L\times Z \rightarrow Z$ by group
homomorphisms, an a cocycle 
$f\in C^2(L,Z,\delta)$, where the differential
$\delta$  is defined by the formula
(\ref{differential1nest}) with repetitions, defines a loop 
satisfying the identity (\ref{reducedcanonical1nest}).
\end{lemma}

First, we compute the product
\begin{align}\label{left_product}
\phi_1\cdots\phi_{n-1}((a_1,z_1),\dots,(a_n,z_n)).
\end{align}Take the product 
$(a_{n-1},z_{n-1})\cdot(a_{n},z_{n})=
(a_{n-1}a_{n},f(a_{n-1},a_{n})+ z_{n-1}a_{n}+z_{n})$
and we multiply on the left:
\begin{align*}
(a_{n-2},z_{n-2})&\cdot(a_{n-1}a_{n},f(a_{n-1},a_{n})+ z_{n-1}a_{n}+z_{n})=
\\
&(\phi_1\phi_2(a_{n-2},a_{n-1},a_{n}),
f(a_{n-2},a_{n-1}a_{n})+z_{n-2}(a_{n-1}a_{n})
\\
&+f(a_{n-1},a_{n}) + z_{n-1}a_{n}+z_{n})=
\\
=&(\phi_1\phi_2(a_{n-2},a_{n-1},a_{n}),
f(a_{n-2},a_{n-1}a_{n})+f(a_{n-1},a_{n})+
\\
&+z_{n-2}(a_{n-1}a_{n}) + z_{n-1}a_{n}+z_{n}).
\end{align*}Further, 
\begin{align}\label{leftproduct_cartesian}
\begin{array}{c}
\phi_1\cdots\phi_{n-1}((a_1,z_1),\dots,(a_n,z_n))=(\phi_1\cdots\phi_{n-1}(a_1,\dots , a_n),\\
\\
\sum_{i=1}^{n-2}f(a_{n-i-1},\phi_1\cdots\phi_{i-1}(a_{n-i},\dots , a_n))+
\sum_{i=1}^{n}z_{i}a_{i+1}\cdots a_{n+1})
\end{array}
\end{align}where $a_{n+1}=e$ the neutral element of the loop $L$.

Now we need to compute the product
\begin{align}\label{right_product}
&\phi_{1}^{u_q}\phi_{2}
\cdots\phi_{t_1-1}\phi_{t_1}^{u_2}\phi_{t_1+1}\cdots 
\phi_{t_1+\cdots +t_{p-1}}^{u_1}
\phi_{m-t_p+1}\cdots\phi_{m-1}((a_1,z_1),\dots,(a_n,z_n)).
\end{align}It is clear from the product formula 
(\ref{productrule1nest}) that the left coordinate of this 
product is
\begin{align*}
&\phi_{1}^{u_q}\phi_{2}
\cdots\phi_{t_1-1}\phi_{t_1}^{u_2}\phi_{t_1+1}\cdots 
\phi_{t_1+\cdots +t_{p-1}}^{u_1}
\phi_{m-t_p+1}\cdots\phi_{m-1}(a_1,\dots, a_n).
\end{align*}So, lets compute the right coordinate.
For this coordinate, the first partial product is 
given by the right coordinate of the formula for pairings 
to the left (\ref{leftproduct_cartesian}). More 
precisely, this product is
\begin{align*}
&
\sum_{i=1}^{t_p-1}f(a_{m-i-1},\phi_1\cdots\phi_{i-1}(a_{m-i},\dots , a_m)+
\sum_{i=m-t_p}^{m-1}z_{i}a_{i+1}\cdots a_{m} + z_m.
\end{align*}Then, when we multiply on the right. The coordinate 
$a_{m+1}$ acts on this product, and we need to add the coordinate
$z_{m+1}$ and the set of factors evaluated in the last product 
at the left argument, and in the new coordinate $a_{m+1}$ at the
right, i.e.
\begin{align*}
&
f(\phi_1\cdots\phi_{t_p-1}(a_{m-t_p},\dots , a_m),a_{m+1})+\\
&
\sum_{i=1}^{t_p-1}f(a_{m-i-1},\phi_1\cdots\phi_{i-1}(a_{m-i},\dots , a_m))a_{m+1}+
\\
&\sum_{i=m-t_p}^{m-1}z_{i}a_{i+1}\cdots a_{m}a_{m+1} + z_ma_{m+1}+
z_{m+1}.
\end{align*}When we multiply by $(a_{m+2},z_{m+2})$ again on the 
right, we obtain
\begin{align*}
&
f(\phi^2_1\cdots\phi_{t_p-1}(a_{m-t_p},\dots , a_{m+1}),a_{m+2})+\\
&
+f(\phi_1\cdots\phi_{t_p-1}(a_{m-t_p},\dots , a_m),a_{m+1})a_{m+2}+\\
&+\sum_{i=1}^{t_p-1}f(a_{m-i-1},\phi_1\cdots\phi_{i-1}(a_{m-i},\dots , a_m))a_{m+1}a_{m+2}+
\\
&+\sum_{i=m-t_p}^{m-1}z_{i}a_{i+1}\cdots a_{m}a_{m+1}a_{m+2} + z_ma_{m+1}a_{m+2}+
z_{m+1}a_{m+2}+z_{m+2}=\\
&\\
&=f(\phi^2_1\cdots\phi_{t_p-1}(a_{m-t_p},\dots , a_m),a_{m+2})+\\
&
+f(\phi_1\cdots\phi_{t_p-1}(a_{m-t_p},\dots , a_m),a_{m+1})a_{m+2}+\\
&+\sum_{i=1}^{t_p-1}f(a_{m-i-1},\phi_1\cdots\phi_{i-1}(a_{m-i},\dots , a_m))a_{m+1}a_{m+2}+
\\
&+\sum_{i=m-t_p}^{m+1}z_{i}a_{i+1}\cdots a_{m+2}+z_{m+2}.
\end{align*}At the end of this series of multiplications on 
the right, we have the term
\begin{align*}
&
f(\phi^{u_1}_1\cdots\phi_{t_p-1}(a_{m-t_p},\dots , a_{m+u_1-1}),a_{m+u_1})+\\
&
f(\phi^{u_1-1}_1\cdots\phi_{t_p-1}(a_{m-t_p},\dots , a_{m+u_1-2}),a_{m+u_1-1})a_{m+u_1}
+\cdots\\
&
\cdots
+f(\phi_1\cdots\phi_{t_p-1}(a_{m-t_p},\dots , a_m),a_{m+1})a_{m+2}\cdots a_{m+u_1}+
\\
&+\sum_{i=1}^{t_p-1}f(a_{m-i-1},\phi_1\cdots\phi_{i-1}(a_{m-i},\dots , a_m))a_{m+1}a_{m+2}\cdots a_{m+u_1}+
\\
&+\sum_{i=m-t_p}^{m+u_1-1}z_{i}a_{i+1}\cdots a_{m+u_1}+z_{m+u_1}.
\end{align*}Then, we multiply by the element 
$(a_{m-t_p-1},z_{m-t_p-1})$ on the left. The first 
coordinate on the right side 
$\phi^{u_1+1}_1\cdots\phi_{t_p-1}(a_{m-t_p},\dots , a_m,a_{m+u_1})$
acts on the element $z_{m-t_p-1}$ giving a shift of one on the sum
\begin{align*}
&\sum_{i=m-t_p}^{m+u_1-1}z_{i}a_{i+1}\cdots a_{m+u_1}+z_{m+u_1},
\end{align*}more accurately, the sum on the next term is 
\begin{align*}
&\sum_{i=m-t_p-1}^{m+u_1-1}z_{i}a_{i+1}\cdots a_{m+u_1}+z_{m+u_1},
\end{align*}and we add the set of factors evaluated in the left on
the coordinate $a_{m-t_p-1}$ and in the right on the last product
$\phi^{u_1+1}_1\cdots\phi_{t_p-1}(a_{m-t_p},\dots , a_m,a_{m+u_1})$.
So, the next term is 
\begin{align*}
&
f(a_{m-t_p-1},\phi^{u_1+1}_1\cdots\phi_{t_p-1}(a_{m-t_p},\dots , a_m,a_{m+u_1}))+\\
&
f(\phi^{u_1}_1\cdots\phi_{t_p-1}(a_{m-t_p},\dots , a_{m+u_1-1}),a_{m+u_1})+\\
&
f(\phi^{u_1-1}_1\cdots\phi_{t_p-1}(a_{m-t_p},\dots , a_{m+u_1-2}),a_{m+u_1-1})a_{m+u_1}
+\cdots\\
&
\cdots
+f(\phi_1\cdots\phi_{t_p-1}(a_{m-t_p},\dots , a_m),a_{m+1})a_{m+2}\cdots a_{m+u_1}+
\\
&+\sum_{i=1}^{t_p-1}f(a_{m-i-1},\phi_1\cdots\phi_{i-1}(a_{m-i},\dots , a_m))a_{m+1}a_{m+2}\cdots a_{m+u_1}+
\\
&+\sum_{i=m-t_p-1}^{m+u_1-1}z_{i}a_{i+1}\cdots a_{m+u_1}+z_{m+u_1}.
\end{align*}Once again, a) left multiplications account for left 
arguments with the last previous product on the right, b)
right multiplications account for right arguments with the last
previous product on the left, and c) the action 
$Z\times L \rightarrow Z$ is used to fill any missing arguments
on the right. Therefore, the last term is the sum of the positive 
side of the differential (\ref{differential1nest}) and the term 
\begin{align}\label{residualsum}
&\sum_{i=1}^{n}z_{i}a_{i+1}\cdots a_{n+1},
\end{align}where $a_{n+1}=e\in L$ is the neutral element of this loop.
More precisely, the right coordinate of the product 
(\ref{right_product}) is
\begin{align*}
& f(\phi_{m-t_1-t_2-\cdots-t_p}^{u_q-1}
    \phi_{m-t_1-t_2-\cdots-t_p+1}
    \cdots
  \phi_{m-t_1-t_2}^{u_2}
   \phi_{m-t_1-t_2+1}\cdots
   \\
&   \cdots
   \phi_{m-t_1-1}
   \phi_{m-t_1}^{u_1}\phi_{m-t_1+1}
   \cdots
   \phi_{m-1}\phi_{m}(a_1,\dots ,a_{n-1}),a_{n})+
   \\
&+   f(\phi_{m-t_1-t_2-\cdots-t_p}^{u_q-2}
    \phi_{m-t_1-t_2-\cdots-t_p+1}
    \cdots
   \phi_{m-t_1-t_2}^{u_2}
   \phi_{m-t_1-t_2+1}\cdots
   \\
  & \cdots
   \phi_{m-t_1-1} 
   \phi_{m-t_1}^{u_1}\phi_{m-t_1+1}
   \cdots
   \phi_{m-1}\phi_{m}(a_1,\dots ,a_{n-2}),a_{n-1})
   a_{n}+\cdots
   \\
 &  \cdots +f(\phi_{1}^{u_1-1}\cdots
   \phi_{m-1}(a_{m-t_1+1}, \dots, a_{m+u_1-2}),
   a_{m+u_1-1})a_{m+u_1}\cdots a_{n}+
\\
&\cdots +
f(a_{m-1}, a_{m})a_{m+1}\cdots a_{n}+\\
&+\sum_{i=1}^{n}z_{i}a_{i+1}\cdots a_{n+1},
\end{align*}

Then, when we substitute in the equation 
\begin{align*}
\phi_1 &\cdots \phi_{n-1}\circ r((a_1,z_1),(a_2,z_2),(a_3,z_3))=\\
=\phi_{1}^{u_q}\phi_{2}&
\cdots\phi_{t_1-1}\phi_{t_1}^{u_2}\phi_{t_1+1}\cdots\\
&\cdots 
\phi_{t_1+\cdots +t_{p-1}}^{u_1}
\phi_{m-t_p+1}\cdots\phi_{m-1}\circ r((a_1,z_1),(a_2,z_2),(a_3,z_3))
\end{align*}we see that the $Q$-coordinates coincide because
$Q$ is a loop of the required class, and in the $Z$-coordinates the 
terms (\ref{residualsum}) cancel.$\hfill \square$

From lemmas \ref{cocycle_2} and \ref{cartesian_extension} it is obtained
the following.

\begin{lemma}
A set of factors $f\in C^2(Q,Z)$ define a loop extension 
with given $Q$-action over $Z$ satisfying the identity 
\ref{canonical1nestlaw} if and only if 
$f$ is a cocycle.
\end{lemma}

The proof of the fact that equivalent extensions
define cocycles that differ by a coboundary does not 
depend on the given identity but only on the 
product formula (\ref{multiplication}) as it is shown in
\cite{Johnson_Leedham-Green}, we only write 
it here to make the text more self-contained.

\begin{lemma}
If the cocycles $f,g\in C^2(Q,Z)$ define
equivalent loop extensions, then they differ by a coboundary, i.e.
$f-g=\delta h$ for some function $h\in C^2(Q,Z)$. 
\end{lemma}
Assume that the cocycles $f,g\in C^2(Q,Z)$ define
equivalent loop extensions $L$ and $M$. This means that
there is a commutative diagram
\begin{align}\label{equivalent_extensions}
\xymatrix{   &                  & L \ar[rd] &          \\
0 \ar[r]     & Z \ar[ru]\ar[rd] &                 & Q \ar[r] & 0.\\
             &                  & M \ar[ru]\ar[uu]^\varphi  &
             }
\end{align}
Then, for sections $s: Q \rightarrow L$,
$u: Q \rightarrow M$ and an element $x\in Q$ the elements 
$x_s,\varphi(x_u)\in L$ and differ by an element $h(x)\in Z$, 
i.e.  $\varphi(x_u)=x_sh(x)$. This condition defines a function
$h:Q\rightarrow Z$ such that $h(e)=e$ (because the sections
are normalized). So, $\varphi(x_uy_u)=
\varphi(x_u)\varphi(y_u)=x_sh(x)y_sh(y)=
  (xy)_sf(x,y)T_{y_s}h(x)h(y)$. But
$x_uy_u= (xy)_ug(x,y)$ and $\varphi((xy)_ug(x,y))=$\linebreak
$=\varphi((xy)_u)\varphi(g(x,y))
=\varphi((xy)_u)g(x,y)=(xy)_sh(xy)g(x,y)$, because 
$\varphi(z)=z$ for every element $z\in Z$. Therefore 
$f(x,y)T_{y_s}h(x)h(y)=h(xy)g(x,y)$, which, in additive
notation gives 
$f(x,y)+h(x)y +h(y)=h(xy)+g(x,y)$ or, equivalently
\begin{equation}\label{equivalent_cocycles}
g(x,y)-f(x,y)=h(y)-h(xy)+h(x)y=\delta h(x,y),
\end{equation}i.e. $g-f=\delta h$. $\hfil \square$

\begin{cor}\label{independence_section}
The cohomology class associated to the extension does 
not depend on the particular section $s:Q\rightarrow L$.
\end{cor}
\proof 
The case $M=L$, $\phi=\id$ gives the needed result. $\hfill \square $

Next, we need to show that, if the set of factors 
is obtained from a given extension (using some section) 
then the product rule (\ref{productrule1nest}) defines an
extension equivalent to the initial one. 

\begin{lemma}\label{Cartesian_abstract_equivalence}
Let $L$ be a loop extension, $s:Q \longrightarrow L$ be a 
normalized section and $f:Q\times Q \longrightarrow Z$ the set 
of factors obtained using this section. The multiplication (\ref{productrule1nest}) on the Cartesian product $Q\times Z$
defines a loop extension equivalent to $L$.
\end{lemma}
\proof\,
Define $\varphi:Q\times Z \longrightarrow L $ by 
the rule $\varphi(a,z)=a_sz$. Then 
$\varphi[(a,w)\cdot(b,z)]= \varphi(ab,f(a,b)+wb+z)=
(ab)_sf(a,b)T_{b_s}(w)z=a_szb_sw=\varphi(a,w)\varphi(b,z)]$.
This means that $\varphi$ is a homomorphism. The 
proof of the commutativity of the diagram 
(\ref{equivalent_extensions}) with maps 
$Z \to  Q\times Z,\;z\mapsto (e,z)$, 
$Q\times Z\to Q,\; (a,z)\mapsto a$ is left as an exercise (or conf.
\cite{Johnson_Leedham-Green}).
$\hfill \square$

Lemmas \ref{cocycle_2}--\ref{Cartesian_abstract_equivalence} and 
corollary \ref{independence_section} constitute the proof of the 
following expected result.

\begin{theorem}\label{classification_extensions_1_nest}
The set of equivalent classes
of extensions of the loop $Q$ belonging to the class defined
by the identity (\ref{reducedcanonical1nest3variables}) 
over the commutative group $Z$ with given action 
$Z\times Q \longrightarrow Z$ is in one-to-one 
correspondence with the cohomology group 
$H^2(Q,Z;\delta)$, where $\delta$ is the differential 
(\ref{differential1nest}) with repetitions.
\end{theorem} 

It may be also true the analogue result for
every class of loops defined by one or several identities.
This issue will be addressed in future works.
\section{Some computations and a problem}
In the sequel, all the action are considered trivial.

\begin{lemma}\label{image_diff_trivial_action}
For trivial action 
$\mathbb{Z}/n\times\mathbb{Z}/m \longrightarrow \mathbb{Z}/m$
the image $\im \delta$ of the differential 
$\delta :C^1(\mathbb{Z}/n, \mathbb{Z}/m )\longrightarrow
C^2(\mathbb{Z}/n, \mathbb{Z}/m )$ given by (\ref{image_diff_trivial_action}) has $(n-1)m/(n,m)$ elements.
\end{lemma}
\proof\, In order to compute the number of elements in the image of
the differential $\delta$ in (\ref{differentialdimension1}), 
note that there are $(n-1)m$ normalized functions of 
one variable, i.e. $(n-1)m$ elements in 
$C^1(\mathbb{Z}/n, \mathbb{Z}/m )$. On the other hand, 
the kernel of $\delta$ in this case is exactly the set 
of group homomorphisms 
$h:\mathbb{Z}/n\longrightarrow \mathbb{Z}/m$, and there 
are exactly $(n,m)$ of them. 
Therefore, there are $(n-1)m/(n,m)$ elements
in the image of $\delta$. $\hfill\square$

\subsection{Left Bol metacyclic loop extensions}
Consider extensions 
\begin{equation}
0\rightarrow \mathbb{Z}/m \rightarrow L\rightarrow \mathbb{Z}/n \rightarrow 0
\end{equation}where $L$ is a left Bol loop. 
In this case we have de differential 
(\ref{left-Bol_diff_no_repetitions}) with $w=y$, i.e.
\begin{equation}\label{left-Bol_differential}
(\delta f)(x,y,z)=f(y\cdot xy,z)+f(w,xy)z+f(x,y)z-
f(y,x\cdot yz)-f(x,yz)-f(y,z).
\end{equation}
Assume that the action 
$T:\mathbb{Z}/n\times \mathbb{Z}/m \longrightarrow \mathbb{Z}/m $ 
is trivial. Then (\ref{left-Bol_diff_no_repetitions}) becomes
\begin{equation}\label{left-Bol_diff_trivial_action}
(\delta f)(x,y,z)=f(y\cdot xy,z)+f(y,xy)+f(x,y)-
f(y,x\cdot yz)-f(x,yz)-f(y,z).
\end{equation}Using additive notation and the commutativity of 
the group $\mathbb{Z}/n$ we obtain
\begin{equation}\label{left-Bol_diff_additive_tr_action}
\begin{array}{l}
(\delta f)(x,y,z)=\\
\hspace{2cm} f(x+2y,z)+f(y,x+y)+f(x,y)-\\
\hspace{3cm}-f(y,x+y+z)-f(x,y+z)-f(y,z).
\end{array}
\end{equation}Note that for normalized sets of factors, this expression vanishes when $x=0$ or $y=0$, but not, in general, when $z=0$.

For a left Bol loop extension,
$(\delta f)(x,y,z)=0$ for every $x,y,z\in \mathbb{Z}/n$. In the case
$x=0$, we have 
$$
\begin{array}{l}
(\delta f)(0,y,z)=\\
\hspace{2cm} f(2y,z)+f(y,y)+f(0,y)-\\
\hspace{3cm}-f(y,y+z)-f(0,y+z)-f(y,z)=0,
\end{array}
$$but, as the set of factors are normalized, $f(0,x)=f(x,0)=0$,
then
\begin{equation}
(\delta f)(0,y,z)=f(2y,z)+f(y,y)-f(y,y+z)-f(y,z)=0.
\end{equation}Therefore,
\begin{equation}
f(y,y)=f(y,y+z)+f(y,z)-f(2y,z).
\end{equation}The right side of this equation does not depend on $z$.
For $z=n-y$, we obtain
\begin{equation}
f(y,y)=f(y,0)+f(y,n-y)-f(2y,n-y)=f(y,n-y)-f(2y,n-y),
\end{equation}i.e.
\begin{equation}\label{diagonal_Bol_metacyclic}
f(y,y)=f(y,n-y)-f(2y,n-y).
\end{equation}

\begin{lemma}\label{reduced_diff_Bol-2-3type}
For $n=3,\; m=2$ equation (\ref{diagonal_Bol_metacyclic})
is equivalent to $\partial f=0$.
\end{lemma}
\proof\, Interpret the set of factors $f$ as a $3\times 3$ 
Boolean matrix with entries $f(x,y)$. We consider all the 
different possible values of the diagonal and compute the 
rest of the matrix entries. As the entries are Boolean, 
we can forget about signs. Note that $f(0,0)=0$. So, 
it is enough to compute for $x\neq 0$, $y\neq 0$. 

$i)\; f(y,y)=0,\; y=1,2$ 

In this case 
$$
0=f(1,1)=f(1,2)+f(2,2)=f(1,2)
$$and
$$
0=f(2,2)=f(2,1)+f(1,1)=f(2,1).
$$ This is the zero matrix.

$ii)\; f(y,y)=1,\; y=1,2$ 

In this case 
$$
1=f(1,1)=f(1,2)+1,
$$so $f(1,2)=0$, and
$$
0=f(2,2)=f(2,1)+1
$$so $f(2,1)=0$. This is the matrix
$
\left(\begin{array}{ccc}
0&0&0\\
0&1&0\\
0&0&1\\
\end{array}\right).
$

$iii)\; f(1,1)=1,\; f(2,2)=0$ 

In this case 
$$
1=f(1,1)=f(1,2)+0,
$$so $f(1,2)=1$, and
$$
0=f(2,2)=f(2,1)+1
$$so $f(2,1)=1$. This is the matrix
$
\left(\begin{array}{ccc}
0&0&0\\
0&1&1\\
0&1&0\\
\end{array}\right).
$

$iv)\; f(1,1)=0,\; f(2,2)=1$ 

This is analogous to the previous case and it is obtained
the matrix
$\left(\begin{array}{ccc}
0&0&0\\
0&0&1\\
0&1&1\\
\end{array}\right).$

To conclude the proof, it is enough to show that the differential 
vanishes for the sets of factors defined by the last three matrices.
As the differential (\ref{left-Bol_diff_additive_tr_action}) 
vanishes when $x=0$ or $y=0$, and when $z=0$ this is equation 
(\ref{diagonal_Bol_metacyclic}), we need only to consider the case
when all the coordinates are non-zero.

Indeed, $
(\delta f)(1,1,1)=
f(0,1)+f(1,2)+f(1,1)
+f(1,0)+f(1,2)+f(1,1)=0
$ because all terms appear twice in the sum;

$
(\delta f)(2,1,1)=
f(1,1)+f(1,0)+f(2,1)
+f(1,1)+f(2,2)+f(1,1)=
f(2,1)+f(2,2)+f(1,1)=0$ by equation (\ref{diagonal_Bol_metacyclic})
for $y=2$;

$
(\delta f)(1,2,1)=
f(2,1)+f(2,0)+f(1,2)
+f(2,1)+f(1,0)+f(2,1)=f(1,2)+f(2,1)=0
$ because for all cases $f(1,2)=f(2,1)$;

$
(\delta f)(1,1,2)=
f(0,2)+f(1,2)+f(1,1)
+f(1,1)+f(1,0)+f(1,2)=
f(1,2)+f(1,1)
+f(1,1)+f(1,2)=0
$ because all remaining terms appear twice;

$
(\delta f)(2,2,1)=
f(0,1)+f(2,1)+f(2,2)
+f(2,0)+f(2,0)+f(2,1)=
f(2,1)+f(2,2)+f(2,2)+f(2,1)=0
$ because all remaining terms appear twice;

$
(\delta f)(2,1,2)=
f(1,2)+f(1,0)+f(2,1)
+f(1,2)+f(2,0)+f(1,2)=f(2,1)
+f(1,2)=0
$ because $f(1,2)=f(2,1)$;

$
(\delta f)(1,2,2)=
f(2,2)+f(2,0)+f(1,2)
+f(2,2)+f(1,1)+f(2,2)=
f(1,2)
+f(2,2)+f(1,1)=0$ by equation (\ref{diagonal_Bol_metacyclic})
for $y=1$;

$
(\delta f)(2,2,2)=
f(0,2)+f(2,1)+f(2,2)
+f(2,0)+f(2,1)+f(2,2)=
f(2,1)+f(2,2)+
f(2,1)+f(2,2)=0
$ because all remaining terms appear twice. $\hfill \square$

\begin{cor}\label{kernel_type_2_3}
The kernel $\ker \delta$ has exactly $4$ elements. 
\end{cor}
\proof\, As all the sets of factors satisfying equation (\ref{diagonal_Bol_metacyclic}) where presented in the proof of
lemma \ref{reduced_diff_Bol-2-3type}, this number is at most 4. 
By the same lemma, all those elements are cocycles. 
$\hfill \square$

\begin{theorem}\label{cohomology_2_3_Bol}
The group $H^2_\delta (\mathbb{Z}/n,\mathbb{Z}/m)$ is trivial. 
\end{theorem}
\proof\, By lemma \ref{image_diff_trivial_action} the image $\im \delta$
has $4$ elements, and by corollary \ref{kernel_type_2_3} the kernel
$\ker \delta$ has the same number of elements.
$\hfill \square$

\begin{theorem}
For $n=3,\; m=2$ there is exactly one left Bol 
loop extension up to equivalence.
\end{theorem}
\proof\, This is a consequence of theorems \ref{classification_extensions_1_nest} and \ref{cohomology_2_3_Bol}.
$\hfill\square$

\subsection{Commutative metacyclic loop extensions}
Consider extensions 
\begin{equation}
0\rightarrow \mathbb{Z}/m \rightarrow L\rightarrow \mathbb{Z}/n \rightarrow 0
\end{equation}where $L$ is a commutative loop. 
In this case the action 
$T:\mathbb{Z}/n\times \mathbb{Z}/m \longrightarrow \mathbb{Z}/m $ is trivial.

From the identity 
$$
x\cdot y = y \cdot x 
$$it is obtained the identity
$$
f(x,y)= f(y,x) 
$$which gives the differential 
$\delta : C^2(\mathbb{Z}/n, \mathbb{Z}/m )\longrightarrow
C^2(\mathbb{Z}/n, \mathbb{Z}/m )$ 
$$
\delta f(x,y)=f(x,y)-f(y,x).
$$That is, $f$ is a cocycle if and only if it is a symmetric 
function. Using analogous methods one can prove that the
number of non-equivalent extensions is given 
by the corresponding cohomology group $H^2_\delta(\mathbb{Z}/n, \mathbb{Z}/m)$. We compute the number of elements in this group.

The set of normalized symmetric functions with values in 
$\mathbb{Z}/m$ is in one-to-one correspondence with
the set of symmetric $(n-1)\times (n-1)$ matrices with entries in this group.
Therefore, there are $\frac{n(n-1)}{2}m$ of this functions.

In order to compute the number of elements in the image of
the differential $\delta$ in (\ref{equivalent_cocycles}), 
note that there are $(n-1)m$ normalized functions of 
one variable, i.e. $(n-1)m$ elements in 
$C^1(\mathbb{Z}/n, \mathbb{Z}/m )$. On the other hand, 
the kernel of $\delta$ in this case is exactly the set 
of group homomorphisms 
$h:\mathbb{Z}/n\longrightarrow \mathbb{Z}/m$, and there 
are exactly $(n,m)$ of them (the great common divisor of this numbers). Therefore, there are $(n-1)m/(n,m)$ elements
in the image of $\delta$. 

We conclude that the number of elements in 
$H^2_\delta(\mathbb{Z}/n, \mathbb{Z}/m)$ is 
$$
\frac{n(n-1)m}{2}\Big/ \frac{(n-1)m}{(n,m)}=\frac{n(n,m)}{2}.
$$

\subsection{Metacyclic loop extensions with inverse property}

Consider extensions 
\begin{equation}
0\rightarrow \mathbb{Z}/m \rightarrow L\rightarrow \mathbb{Z}/n \rightarrow 0
\end{equation}where $L$ is a loop with unique inverse property. 
In this case the action 
$T:\mathbb{Z}/n\times \mathbb{Z}/m \longrightarrow \mathbb{Z}/m $ might be non trivial.

Then, the multiplication rule in 
$L\approx \mathbb{Z}/n\times \mathbb{Z}/m$ is given by
$$
(x,a)\cdot(y,b)=(x+y, f(x,y)+ at^y+b),
$$where $T(1)=t$, the image of the generator under $T$, and the action is by integer multiplication.

Because the loop $L$ has unique inverses, for a given pair
$(x,a)$, we have an element $(y,b)$ such that
$$(x+y, f(x,y)+ at^y+b)=(0,0)$$
and
$$(y+x, f(y,x)+ bt^x+a)=(0,0).$$
So, $y=-x$. Then, $f(x,-x)+ at^{-x}+b=0$ and 
$f(-x,x)+ bt^{x}+a=0$. Multiplying in the last
equation by $t^{-x}$, we obtain 
$f(-x,x)t^{-x}+ b+at^{-x}=0$. Therefore,
$$
f(x,-x)-f(-x,x)t^{-x}=0.
$$

Define the operator $\delta : 
C^2(\mathbb{Z}/n, \mathbb{Z}/m )\longrightarrow C^2(\mathbb{Z}/n, \mathbb{Z}/m )$ by 
$$ \delta f (x,y)=f(x,y)-f(y,x)y
$$and consider the map $\phi: \mathbb{Z}/n \longrightarrow
\mathbb{Z}/n\times \mathbb{Z}/n; x \mapsto (x,-x)$.

Then, the cocycle condition for this kind of loop is
$(\delta f) \circ \phi =0$. Lets compute $\delta^2$ in 
this case. Take 
$f\in C^1(\mathbb{Z}/n, \mathbb{Z}/m)$, set
$h(x,y)=\delta f(x,y)= f(y)-f(xy)+f(x)y$. Then
\begin{align*}
\delta\delta f (x,y)=\delta h (x,y)&=h(x,y)-h(y,x)y\\
&= f(y)-f(xy)+f(x)y-(f(x)-f(yx)+f(y)x)y=\\
&=f(y)-f(xy)+f(x)y-f(x)y+f(yx)y-f(y)xy=\\
&=f(y)-f(xy)+f(yx)y-f(y)xy.
\end{align*}This does not vanish in general. Now 
substitute $y=x^{-1}$ (we use multiplicative notation
here, as this is valid for every loop extension of 
the given type) to obtain
\begin{align*}
(\delta\delta f) \circ \phi (x)&=
\delta\delta f (x,x^{-1})=\\
&=f(x^{-1})-f(xx^{-1})+f(x^{-1}x)x^{-1}-f(x^{-1})xx^{-1}=\\
&=f(x^{-1})-f(e)+f(e)x^{-1}-f(x^{-1})e=\\
&=f(x^{-1})-f(x^{-1})=0.
\end{align*}

%\subsection*{Trivial action}
%Now we compute the number metacyclic extensions with inverse
%property and trivial action.

%Assuming that the action 
%$T:\mathbb{Z}/n\times \mathbb{Z}/m \longrightarrow \mathbb{Z}/m $ is trivial, we obtain the same differential that we obtained in the commutative case, but it only vanishes when restricted to
%the diagonal. This is not strange, because the unique inverse property means that an element commutes with its inverse.

%In this case, as the square of the differential does not vanish, 
%we can not proceed as before to count the number of extensions.

So, take the set 
$U= \Delta (\mathbb{Z}/n)\subset \mathbb{Z}/n\times \mathbb{Z}/n$ and let 
$r: C^2(\mathbb{Z}/n, \mathbb{Z}/m)\longrightarrow C^2(U, \mathbb{Z}/m)$ be the restriction homomorphism, i.e.
$r(f)=f\vert_U$. 

%Then the following diagram

%\begin{equation}
%\xymatrix{
%C^{1}\ar[r]^{\delta}\ar[d]_{||} &\ar[r]^{\delta}C^{2}\ar[d]_{r} %&C^2\ar[d]  \\
%C^{1}\ar[r]^{\delta}            &C^{1}\ar[r]^{\delta}& C^1  \\
% }\end{equation}
%is commutative. As in the upper row $\delta^2\neq 0$, this row does not %define
%any homology. 

Consider the commutative diagram
\begin{equation}
\xymatrix{
C^{1}\ar[r]^{\delta} &\ar[r]^{\delta}C^{2}\ar[d]_{r} &0  \\
        &C^{1}& \\
 }\end{equation}The upper row, according to \cite{Johnson_Leedham-Green} computes the number of arbitrary loop extensions. So, in order to compute the number of extensions with unique inverse property, one has to find
all the set of factors $f\in C^2(\mathbb{Z}/n, \mathbb{Z}/m)$ such 
that $(\delta\delta f) \circ \phi (x)=0$ and then compute the number of
corresponding classes $\sigma_f\in L^2(\mathbb{Z}/n, \mathbb{Z}/m)$.
Does not seems to be any clear method to compute this number.


\begin{thebibliography}{20}
\bibitem{Eilenberg-Maclane} Eilenberg S., Maclane S., \textit{Cohomology 
theory in abstract groups I} 
Ann. of math. Vol. 48, No. 1 (1947), 51--78.
\bibitem{Nishigori_1} Nishigori N., 
\textit{On Loop Extensions of Groups and M-cohomology Groups I.} 
J. Sci. HIROSHIMA UNIV. SER. A-I 27 (1963), 151--165.\bibitem{Nishigori_2} Nishigori N., 
\textit{On Loop Extensions of Groups and M-cohomology Groups II.} 
J. Sci. HIROSHIMA UNIV. SER. A-I 29 (1965), 17--26.
\bibitem{Johnson_Leedham-Green} Kenneth W.J., Ch.R. Leedham-Green
\textit{Loop cohomology} Czec. Math. J., Vol. 40 (1990), No. 2, 182--194.
\bibitem{Nagy_Strambach} P.T. Nagy; K. Strambach
\textit{Schreier loops}
Czec. Math. J., Vol. 58 (2008), No. 3, 759--786.
\end{thebibliography}
\end{document}